\documentclass[12pt]{article}

\usepackage[latin1]{inputenc}
\usepackage{color}
\usepackage{latexsym}
\usepackage{amssymb}
\usepackage{graphicx}
\usepackage{enumerate}
\usepackage{extarrows}
\newtheorem{Theorem}{Theorem}[part]
\newtheorem{Definition}{Definition}[part]
\newtheorem{Proposition}{Proposition}[part]

\newtheorem{Corollary}{Corollary}[part]
\newtheorem{Remark}{Remark}[part]

\def \s{~~~}
\def \2{\vspace{2mm}}

\def \R{\mathbb{R}}

\def \E{\mathbb{E}}
\def \F{\mathbb{F}}

\def \K{\tilde{K}}
\def \Y{\tilde{Y}}
\def \Kh{\hat{K}}

\def \mA{\mathcal{A}}

\def \Dc{{\cal D}}

\def \Fc{{\cal F}}

\def \ep{\hbox{ }\hfill$\Box$}

\def\beqs{\begin{eqnarray*}}
\def\enqs{\end{eqnarray*}}
\def\beq{\begin{eqnarray}}
\def\enq{\end{eqnarray}}

\begin{document}

\title{Skorohod equation and BSDE's with two reflecting barriers}
\author{Soufiane Aazizi \thanks{  Department of Mathematics, Faculty of Sciences Semlalia Cadi Ayyad University, B.P. 2390 Marrakesh, Morocco. \texttt{aazizi.soufiane@gmail.com}} }

%\date{}

\maketitle

\begin{abstract}
We solve a class of doubly reflected backward stochastic differential equation whose generator depends on the resistance due to reflections, which extend the recent work of Qian and Xu on reflected BSDE with one barrier. We then obtain the existence and uniqueness and the continuous dependence theorem for this reflected BSDE.
\end{abstract}

\vspace{13mm}

\noindent {\bf Key words~:} Local Time; Reflected backward stochastic differential equation; Picard iteration; Optional projection; Extended Skorohod problem.

\vspace{5mm}


\newpage
\section {Introduction}
In this paper we are interested to study the following doubly reflected backward stochastic differential equation (DBSDE for short)
\begin{equation}
\label{RBSDE}
\left\{
\begin{array}{lll}
Y_t = \xi  +\int_t^Tf(r, Y_r, Z_r,K_r )dr - \int_t^T Z_rdW_r + K_T - K_t\\
K=K^l-K^u \\
\forall t \leq T, \s L_t\leq Y_t \leq U_t \s \mbox{and} \s \int_0^T(U_r - Y_r)dK^u_r=\int_0^T(Y_r- L_r)dK^l_r=0
\end{array}
\right.
\end{equation}
where $W$ is a Brownian motion defined on some complete filtered probability space $(\Omega, \mathcal{F}, P, )$, a
terminal value $\xi \in \mathcal{F}_T$, a continuous barriers $L$ and $U$ which are modeled by a semimartingales, and $K^l$ and $K^u$ are nondecreasing
processes.
\2 \\
The BSDE was firstly initiated by Bismut \cite{B76} and later developed by Pardoux and Peng \cite{PP90} to prove existence and uniqueness of adapted solution,
 EL Karoui et al. \cite{EKPPQ97} introduced the notion of reflected BSDEs with lower barrier, in which the component Y is forced to stay above a given obstacle,
  Cvitanic, Karatzas and Soner \cite{CKS98}, and later Aazizi and Ouknine \cite{AO11} considered the case where the constraint is imposed on the component Z. In
  the same frame, the generalization of BSDE with two continuous reflecting barriers is introduced by Cvitanic-Karatzas \cite{CK96}. Since then, there were many
  works on the latter kind of BSDEs.
A new kind of reflected BSDE has been introduced by Bank and El Karoui \cite{BK04}, by a variation of Skorohod's obstacle problem, known as variant reflected BSDE, which takes the following forme:
\begin{equation}
\left\{
\begin{array}{lll}
Y_t = X_T  +\int_t^Tf(r, Y_r, Z_r,A_r )dr - \int_t^T Z_rdW_r \\
Y \leq X \\
 \int_t^T|X_r - Y_r|dA_r=0
\end{array}
\right.
\end{equation}
where $A$ is an increasing process, with $A_{0-} = -\infty$. The process $A$ does not directly act on $Y$ to push the solution downwards such that $Y_t \leq X_t$ like in standard BSDE with one barrier, but it acts through the generator $f$. This work has been generalized by Ma and Wang in \cite{MW09}, to prove that a solution in a small-time duration, under some extra conditions, exists and is unique. Recently in the same framework, Qian and Xu \cite{QX11} studied the following  class of reflected BSDE
\begin{equation}
\left\{
\begin{array}{lll}
Y_t = \xi_T  +\int_t^Tf(r, Y_r, Z_r,K_r )dr - \int_t^T Z_rdW_r + K_T - K_t \\
S_t \leq Y_t \\
 \int_t^T(S_r - Y_r)dK_r=0
\end{array}
\right.
\end{equation}
Here, $K$ appears in the driver as a resistance force, if $f$ is decreasing in $K$, then we get an extra force from the Lebesgue integral, if $f$ is increasing in $K$, then there is a kind of cancelation of the positive force, in general case, they consider this RBSDE as an equation with resistance. They derived an explicit formula of the increasing process $K$ by using the result on Skorohod equation together with the theory of optional dual projection (see \cite{HWY92}). \2 \\
In this paper, we extend the approach of \cite{QX11}, to the  doubly reflected BSDE (\ref{RBSDE}). The paper is organized as follows,In Section 2, we provide and explicit formula of the process $K:= K^l - K^u $  using Tanaka formula (see \ref{TanakaK}) and  the extended Skorohod map (ESM in short) (see (\ref{Kt})) introduced  by Ramanan \cite{R06}.  In Section 3 we study a doubly reflected BSDE with generator depending on $K$, we then prove existence and uniqueness through fixed point theorem. Finally, in Section 4, we obtain the continuous dependance property of the solution.

\section{Explicit formula of the process K}

\setcounter{equation}{0} \setcounter{Assumption}{0}
\setcounter{Theorem}{0} \setcounter{Proposition}{0}
\setcounter{Corollary}{0} \setcounter{Lemma}{0}
\setcounter{Definition}{0} \setcounter{Remark}{0}

\subsection{General formulation }

On a given complete probability space $(\Omega, \mathcal{F}, P, )$, let $W:=(W^1, W^2, ..., W^d)$ be a d-dimensional Brownian motion defined on a finite interval $[0, T]$, and denote  $\F= \{\mathcal{F}_t \}_{0 \leq t \leq T}$ the filtration  generated by the Brownian motion and the $P$-null sets. We denote by $\mathcal{P}$ the progressive $\sigma$-field on the product space $[0, T] \times \Omega $. \2 \\
We consider the following spaces:
\begin{itemize}
  \item $\mathcal{L}^2 (\mathcal{F}_t)$ the space of all $\mathcal{F}_t$-measurable real random variable $\phi$ such that $$\mathbb{E}|\phi|^2 < \infty.$$
  \item $H^2_d(0, T)$ the space of $\mathbb{R}^d$-valued predictable process $\psi$ such that
  $$\mathbb{E} \left[ \int_0^T|\psi_t|^2dt\right] < \infty.$$
  \item $\mathcal{S}^2(0, T)$ is the space of all continuous semimartingales over $(\Omega, \mathcal{F}, P,\{\mathcal{F}_t \}_{0 \leq t \leq T} )$.
  \item $\mathcal{A}^2(0, T)$ the space of all $\mathcal{F}_T$-measurable continuous and increasing process $K$ with $K_0 =0$ and such that $\mathbb{E}|K_T|^2 < \infty$.
  \item $D_\mathcal{F}^2(0,T)$ the set of $\mathbb{F}$-progressively measurable càdlàg process $\phi$ with     $$ \mathbb{E}\left[ \sup_{0 \le t \leq T} |\phi_t|^2 \right] < \infty.$$
\end{itemize}
In the sequel of this paper, we denote
\begin{enumerate}[(i)]
  \item The process $K$ to be the difference of the two increasing processes $K^l$ ans $K^u$ such that $K:= K^l - K^u$ with $K^l$, $K^u$ $\in \mA^2$.
  \item The generator $f:[0, T]\times \Omega \times S^2 \times H^2 \times \mA^2\times \mA^2 \rightarrow \R$ is globally Lipschitz w.r.t $y,z$ and $k$ such that:
        \begin{eqnarray}\label{f}
        |f(r,y,z,k) - f(r,y',z',k')| &\leq &L_1 \big(|y- y'| + |z-z'|\big) + L_2|k-k'|.
        \end{eqnarray}
\end{enumerate}
We recall the definition of optional (dual) projection (See Nikeghbali \cite{N06} for more details).
\begin{Definition} \label{DefOP}{\emph{Optional projection}} \\
Let $X$ be a measurable process either positive or bounded. There exist a unique (up to indistinguishability) optional process $X^\flat$ called \emph{optional projection} such that:
\begin{equation}
\mathbb{E}[X_T 1_{T<\infty}/\mathcal{F}_T ] = X^\flat_T 1_{T< \infty} \s \mbox{ a.s.}
\end{equation}
for every stopping time T.
\end{Definition}
\begin{Definition} \label{DefDOP}{\emph{Dual optional projection}} \\
Let $(A_t)$ be an integrable raw increasing process. We call \emph{dual optional projection} of $A$  the (optional) increasing process $(A^o_t)$ defined by:
\begin{equation}
\mathbb{E} \left[ \int_{\mathbb{R}^+} X_r dA^o_r\right]= \mathbb{E} \left[ \int_{\mathbb{R}^+} X^\flat_r dA_r\right],
\end{equation}
for any bounded measurable $X$.
\end{Definition}
By the definitions above, we define the optional projection of $K^l$ (resp. $K^u$) by ${K^l}^\flat$ (resp. ${K^u}^\flat$) and its dual optional projection by ${K^l}^o$
(resp. ${K^u}^o$) which is continuous and increasing. Furthermore, we have ${K^l}^\flat - {K^l}^o$(reps. ${K^u}^\flat - {K^u}^o$) a continuous martingale.
\2 \\
In the first time, we propose to derive some explicit formulas of the process $K$, when the doubly reflected BSDE has a generator $f$ not depending on $y$, $z$ and $k$ taking the following form
\begin{equation} \label{DRBSDEF}
\left\{
\begin{array}{lll}
Y \in S^2, Z\in \mathcal{H}^2_d \mbox{ and } K^l, K^u \in \mA^2 \\
Y_t = \xi  +\int_t^Tf_rdr - \int_t^T Z_rdW_r + K_T - K_t\\
L_t\leq Y_t \leq U_t, \s  \forall t \leq T \\
\mbox{The Skorohod conditions hold: }\int_0^T(U_r - Y_r)dK^u_r=\int_0^T(Y_r- L_r)dK^l_r=0.
\end{array}
\right.
\end{equation}
where $(f_t)_{0 \leq t \leq T}$ is optional, $\mathbb{E} \Big(\int_0^T f^2_r dr\Big)$  $<$  $\infty$, and $L$, $U$ are two continuous semimartingale.
\\
Observe that $Y$ could be written as
\beq \label{RBSDEF}
Y_t = Y_0 - \int_0^tf_rdr + \int_0^tZ_rdW_r- K^l_t + K^u_t.
\enq
In this whole paper, we suppose that
\beq \inf_{t\geq 0} \big(U_t - L_t)>0 \s \mbox{ a.s.}\enq

\subsection{Local time}

The main objective in this part is to provide an explicit formula of the two increasing processes $K^u$ and $K^l$ in term of local time associated to $Y$, based on Tanaka formula. According to \cite{QX11}, we can interpret the increasing processes $K^l$ and $K^u$ in this framework, as a time-reversed local time, in order that $K^l$ (resp. $K^u$) will be called the reflected local time of $Y$ at $L$ (resp. at $U$).

\begin{Corollary}
Assume that $L\leq Y \leq U$ are three continuous semimartingales, such that $L$ and $U$ take the form $X = M^X + A^X $ with $X:= L, U$ and   $(M^X, A^X) \in \mathcal{M}^2 \times  BV[0, \infty) $, when $Y$ is solution to (\ref{RBSDEF}). Then
\begin{eqnarray}\label{TanakaK}
K_t =&-&\int_0^t\left( I_{Y_r=L_r} + I_{Y_r=U_r}\right)f_rdr \nonumber \\
    &-&\int_0^tI_{Y_r=U_r}dA^U_r - \int_0^tI_{Y_r=L_r}dA^L_r + L_t^{U-Y} - L_t^{Y-L},
\end{eqnarray}
with $L^X$ denote the local time of the continuous semimartingale $X$ at $0$.
\end{Corollary}
{\bf Proof. }
From equation (\ref{RBSDEF}), and the following  Tanaka formula
\begin{equation}
\left\{
\begin{array}{lll}
(Y_t-L_t)^-= (Y_0-L_0)^- - \int_0^tI_{Y_r \leq L_r}d(Y_r-L_r) + L_t^{Y-L}\\
(U_t-Y_t)^-=(U_0-Y_0)^- - \int_0^tI_{U_r \leq Y_r}d(U_r-Y_r) + L_t^{U-Y}
\end{array}
\right.
\end{equation}
Together with the fact that $ (Y-L)^-=(U-Y)^-=0$ and increasing property of the local time, we have:
\begin{equation}
\left\{
\begin{array}{lll}
L_t^{Y-L} = - \int_0^tI_{Y_r = L_r}f_rdr - \int_0^tI_{Y_r = L_r}dK_r - \int_0^tI_{Y_r = L_r}dA^L_r \\
L_t^{U-Y} = \int_0^tI_{Y_r = U_r}f_rdr - \int_0^tI_{Y_r = U_r}dK_r + \int_0^tI_{Y_r = U_r}dA^U_r
\end{array}
\right.
\end{equation}
Then
\begin{equation}
\left\{
\begin{array}{lll}
K^l_t = - \int_0^tI_{Y_r = L_r}f_rdr  - \int_0^tI_{Y_r = L_r}dA^L_r - L_t^{Y-L}\\
K^u_t = \int_0^tI_{Y_r = U_r}f_rdr  + \int_0^tI_{Y_r = U_r}dA^U_r - L_t^{U-Y}
\end{array}
\right.
\end{equation}
Since $K = K^l - K^u$, we come at the end of the proof.
\ep \\

\subsection{Extended Skorohod problem}

We derive an explicit formula of the process $K$ solution of doubly reflected BSDE (\ref{RBSDEF}), using Skorohod equation. Throughout this part, $D[0, \infty)$ will denote real-valued càdlàg functions on $[0, \infty)$, $D^-[0, \infty)$ (resp. $D^+[0, \infty)$) will denote càdlàg functions on $[0, \infty)$ taking values in $\mathbb{R}\cup \{-\infty\}$ (reps. in $\mathbb{R}\cup \{\infty\}$) and   $BV[0, \infty)$  denotes the subspace of functions with bounded variation on every finite interval. According to Skorohod equation with two time-boundaries, we have the following definition.

\begin{Definition}(Skorohod problem)\\
Let $\alpha, \beta \in D[0, \infty)$ be such that $\alpha \leq \beta$. Given $x \in D[0, \infty)$ a pair of functions $(y,\eta) \in D[0, \infty) \times BV[0, \infty)$ is said to be a solution of the Skorohod problem on $[\alpha, \beta]$ for $x$ if the following two properties are satisfied:
\begin{enumerate}
  \item $y_t=x_t + \eta_t \in [\alpha_t, \beta_t]$, $\forall t \geq 0$.
  \item $\eta(0^-)=0$, and $\eta$ have the decomposition $\eta:=\eta^l-\eta^u$, where $\eta^l,\eta^u \in I[0, \infty)$,
  \begin{eqnarray}
  \int_0^\infty I_{y_s<\beta_s }d\eta^u_s=\int_0^\infty I_{y_s> \alpha_s}d\eta^l_s=0.
  \end{eqnarray}
\end{enumerate}
\end{Definition} \vspace{2mm}
A more general Skorohod problem is called Extended Skorohod Ptoblem (ESP in short) firstly introduced by Ramanan \cite{R06} (see Definition 2.2 in \cite{S10}) which allows a pathwise construction of reflected Brownian motion that is not necessarily semimartingales.  More recently, Burdzy et al. have shown in Theorem 5 in \cite{BKR09}, that for any $\alpha \in D^-[0, \infty)$ and $\beta \in D^+ [0, \infty)$ such that $\alpha\leq \beta$, there is a well defined Extended Skorohod Map (ESM in short) $\Gamma_{\alpha, \beta}$ represented by
\begin{equation} \label{SkorEq}
\Gamma_{\alpha, \beta}(x_t)=x_t - \Xi_{\alpha, \beta}(x_t)n
\end{equation}
where $\Xi_{\alpha, \beta}(x):D[0, \infty) \rightarrow D[0, \infty) $ is given by
\begin{eqnarray}
\Xi_{\alpha, \beta}(x_t)&=&\max\left\{\left[(x_0 - \beta_0)^+ \wedge \inf_{0 \leq r \leq t}(x_r - \alpha_r)\right];\right. \nonumber\\
                                && \s\s \left.\sup_{0\leq s \leq t}\left[(x_s - \beta_s)^+ \wedge \inf_{s \leq r \leq t}(x_r -
                                \alpha_r)\right]\right\}.
\end{eqnarray}
Note that if $(y, \eta) $ is a solution of the Skorohod problem (SP) on $[\alpha,\beta]$ for $x$, then it is also a solution of extended Skorohod problem (ESP) on $[\alpha, \beta] $ for $x$. However, the expression of the process $\Xi$ is slightly complicate to handle.
\2 \\
In his paper, Slaby \cite{S10} provided an alternative formula (see (\ref{Xi})) for the two sided Skorohod map with time depended boundaries, that is easier to understand and has more interesting properties, especially the Lipschitz property of $\Gamma_{\alpha, \beta}$. Those results are reminded below.
\2 \\
Let us introduce the following notations:
\begin{itemize}
  \item Two pairs of times depending on $x$:
        \begin{eqnarray}
        T_\alpha:= \min\{t>0, / \alpha_t - x_t \geq 0\}, \nonumber \\
        T^\beta:= \min\{t>0, / x_t -\beta_t \geq 0\}. \label{T}
        \end{eqnarray}
  \item Two functions:
    \begin{eqnarray}
    H_{\alpha, \beta}(x_t)&:=& \sup_{0 \leq s \leq t}  \left\{(x_s - \beta_s) \wedge \inf_{s \leq r \leq t}(x_r - \alpha_r) \right\}, \nonumber\\
    J_{\alpha, \beta}(x_t)&:=&\inf_{0\leq s \leq t}\left\{(x_s - \alpha_s) \vee \sup_{s \leq r \leq t}(x_r -\beta_r)\right\}.
    \end{eqnarray}
\end{itemize}
\begin{Corollary}\label{SlabyExpressionXi}
Let $\alpha$ $\in$ $D^-[0, \infty)$ and $\beta$ $\in$ $D^+[0, \infty)$ be such that $\inf_{t \geq 0}(\beta_t - \alpha_t) >0$. Then for every $ x\in D[0, \infty)$, we have:
    \begin{eqnarray} \label{Xi}
        \Xi_{\alpha, \beta}(x_t)=I_{\{T^\beta < T_\alpha\}}I_{[T^\beta, \infty)}(t)H_{\alpha, \beta}(x_t) + I_{\{T_\alpha < T^\beta\}}I_{[T_\alpha, \infty)}(t)J_{\alpha, \beta}(x_t).
    \end{eqnarray}
\end{Corollary}
\2
\begin{Theorem}\label{SlabyLipsh}\emph{Lipschitz continuity}  \\
Under the same conditions of corollary \ref{SlabyExpressionXi}, we have for any $x, x' \in D[0, \infty)$
\begin{equation}
\|\Gamma_{\alpha, \beta}(x) - \Gamma_{\alpha, \beta}(x')  \| \leq  \|x-x'  \|,
\end{equation}
where $\|x \|= \sup_{0 \leq t \leq T} |x_t| $.
\end{Theorem}
\2
According to the doubly reflected BSDE (\ref{DRBSDEF}) and the expression (\ref{RBSDEF}) of $Y$, denote
\beqs \Gamma_{L,U}(x_t) := Y_{T-t}, \s \s \Xi_{L,U}(x_t)=K_{T-t} - K_T,\enqs
where
$x_t = \xi  +\int_{T-t}^Tf_rdr - \int_{T-t}^T Z_rdW_r.
$ \2 \\
It follows from (\ref{SkorEq}) that $K_t = \Xi_{L,U}(x_{T-t}) - \Xi_{L,U}(x_T) $ and more explicitly, we have
\begin{eqnarray} \label{Kt}
K_t &=&  I_{\{T^U < T_L\}}I_{[T^U, \infty)}(t)\left[H_{L,U}(x_{T-t})-H_{L,U}(x_T) \right] \nonumber \\
    && + I_{\{T_L < T^U\}}I_{[T_L, \infty)}(t)\left[J_{L,U}(x_{T-t})-J_{L,U}(x_T) \right].
\end{eqnarray}

\begin{Proposition}\label{PropK}
The two continuous processes $K^u$ and $K^l$ are increasing. Moreover, the measures $dK^u$ and $dK^l$ are carried by the sets $\{Y=U\}$ and $\{Y=L\}$ respectively.
\end{Proposition}
{\bf Proof. }
We can write the process $K^u$ as:
\beqs
K^u_t=\int_0^tI_{\{Y_r=U_r\}}dK_r.
\enqs
According to Skorohod condition in (\ref{DRBSDEF}) we have
\beqs \int_0^T(U_r - Y_r)dK^u_r= \int_0^T(U_r - Y_r)I_{\{Y_r=U_r\}}dK_r=0\enqs
 which means that the support of the measure $dK^u$ is carried by the set $\{Y=U\}$. Similarly, the measure $dK^l$ associated to the process
 $K^l_t=\int_0^tI_{\{Y_r=L_r\}}dK_r$ is carried by the set $\{Y=L\}$.   \\
 From other side, if $T_L < T^U $, then by (\ref{T}) and Skorohod condition in (\ref{DRBSDEF}) is the process
 $K^l$ who plays a role in making $Y$ above the barrier $L$, and according to the explicit formula (\ref{Kt}) of the process $K=K^l-K^u$, we observe that
\begin{eqnarray}
        K^l_t &=&  I_{[T_L, \infty)}(t)\left[\inf_{0\leq s \leq T-t}\left\{(x_s - L_s) \vee \sup_{s \leq r \leq T-t}(x_r -U_r)\right\} \right. \nonumber \\
        &&\left. \ \ \ \ \ \ \ \ \ \ \ \ \ \ \ \ \ \ \ \ \  - \inf_{0\leq s \leq T}\left\{(x_s - L_s) \vee \sup_{s \leq r \leq T}(x_r -U_r)\right\}\right],
    \end{eqnarray}
which is an increasing process. \2 \\
If $T^U < T_L$, then $K^u$ plays the role in making $Y$ below the barrier $U$ such that :
\begin{eqnarray}
        K^u_t &=&  -I_{[T^U, \infty)}(t)\left[\sup_{0 \leq s \leq T-t}  \left\{(x_s - U_s) \wedge \inf_{s \leq r \leq T-t}(x_r - L_r) \right\}\right. \nonumber \\
                                && \left.\ \ \ \ \ \ \ \ \ \ \ \ \ \ \ \ \ \ \ \ \  - \sup_{0 \leq s \leq T}  \left\{(x_s - U_s) \wedge \inf_{s \leq r \leq T}(x_r - L_r) \right\}  \right], \nonumber \\
    \end{eqnarray}
which is  increasing.
\ep \\

%%%%%%%%%%%%%%%%%%%%%%%%%%%%%%%%%%%%%%%%%%%%%%%%%%%%%%%%%%%%%%%%%%%%%%%%%%%%%%%%%%%%%%%%%%%%%
%%%%%%%%%Existence of solution of the RBSDEs with resistance by Picard iteration %%%%%%%%%%%%
%%%%%%%%%%%%%%%%%%%%%%%%%%%%%%%%%%%%%%%%%%%%%%%%%%%%%%%%%%%%%%%%%%%%%%%%%%%%%%%%%%%%%%%%%%%%%

\section{Doubly reflected BSDEs with resistance}
\setcounter{equation}{0} \setcounter{Assumption}{0}
\setcounter{Theorem}{0} \setcounter{Proposition}{0}
\setcounter{Corollary}{0} \setcounter{Lemma}{0}
\setcounter{Definition}{0} \setcounter{Remark}{0}
In this section, we prove existence and uniqueness of a class of doubly reflected BSDE with resistance, by constructing a Picard iteration. We formulate this class of BSDE as the following.

\begin{Definition}\label{DefDRBSDER}
A solution of BSDE with resistance reflected between lower barrier $L \in S^2$ and upper barrier $U \in S^2$ associated to $(\xi, f)$ is a
quadruple $(Y, Z, K^l, K^u)$ $\in$ $\mathcal{D}:= \mathcal{S}^2 \times H^2_d \times \mathcal{A}^2 \times \mathcal{A}^2$ satisfying

\begin{enumerate}[(i)]
  \item $(Y, Z, K^l, K^u)$ solves the following BSDE on $[0, T]$:
        \beq \label{DRBSDER}Y_t = \xi  +\int_t^Tf(r, Y_r, Z_r,K^l_r - K^u_r)dr - \int_t^T Z_rdW_r + K^l_T - K^l_t-(K^u_T - K^u_t). \enq
  \item  $ L_t\leq Y_t \leq U_t, \mbox{ a.s. a.e }\s  \forall t \leq T.$
  \item  Skorohod conditions hold: \beq \label{Skorohod}
        \int_0^T(U_r - Y_r)dK^u_r=\int_0^T(Y_r- L_r)dK^l_r=0, \s \mbox{ a.s.}
        \enq
\end{enumerate}
\end{Definition}
According to (\ref{Kt}), if $(Y,Z,K^l,K^u) $ is solution to reflected BSDE with resistence in the sense of Definition \ref{DefDRBSDER}, then $K_t$ must be:
    \begin{eqnarray} \label{KtSlaby}
        K_t &=&  I_{\{T^U < T_L\}}I_{[T^U, \infty)}(t)\left[\sup_{0 \leq s \leq T-t}  \left\{(x_s - U_s) \wedge \inf_{s \leq r \leq T-t}(x_r - L_r) \right\}\right. \nonumber \\
                                && \left.\ \ \ \ \ \ \ \ \ \ \ \ \ \ \ \ \ \ \ \ \ \ \ \ \ - \sup_{0 \leq s \leq T}  \left\{(x_s - U_s) \wedge \inf_{s \leq r \leq T}(x_r - L_r) \right\}  \right] \nonumber \\
        &+& I_{\{T_L < T^U\}}I_{[T_L, \infty)}(t)\left[\inf_{0\leq s \leq T-t}\left\{(x_s - L_s) \vee \sup_{s \leq r \leq T-t}(x_r -U_r)\right\} \right. \nonumber \\
        &&\left. \ \ \ \ \ \ \ \ \ \ \ \ \ \ \ \ \ \ \ \ \ \ \ \ \  - \inf_{0\leq s \leq T}\left\{(x_s - L_s) \vee \sup_{s \leq r \leq T}(x_r -U_r)\right\}\right],
    \end{eqnarray}
where
\begin{equation}\label{x}
x_t = \xi  +\int_{T-t}^Tf(r, Y_r, Z_r, K^l_r- K^u_r)dr - \int_{T-t}^T Z_rdW_r,
\end{equation}
\subsection{Existence and uniqueness By Picard iteration}

One approach to prove existence of the solution of reflected BSDE's with two barriers, is to use the solution of Skorohod problem by constructing a
Picard-type iterative procedure (See e.g. \cite{EC78} or \cite{EKPPQ97}) to the reflected BSDE with resistance. Throughout this section, we adapt the new method of \cite{QX11} to our setting.
\2 \\
%%%%%%%%% P r o p o s i t i o n  %%%%%%%%%%%%
The following proposition state the mapping which leads to prove existence and uniqueness of the solution $(Y, Z, K^l, K^u)$.
\begin{Proposition} \label{phi} {\textbf{Picard iteration}} \\
The mapping $\phi: \mathcal{D}:= \mathcal{S}^2 \times H^2_d \times \mathcal{A}^2 \times \mathcal{A}^2\rightarrow \mathcal{D}$ which associates $(Y, Z, K^l, K^u)  $ to $\phi(Y, Z, K^l, K^u) =(\tilde{Y}, \tilde{Z}, \tilde{K}^l,\tilde{K}^u ) $ is well defined. Moreover, the decomposition of $\tilde{Y}$ is given by:
\beq \label{phi(Y)} \tilde{Y}_t = \xi + \int_t^T f(r, Y_r, Z_r, {K^l}^\flat_r-{K^u}^\flat_r)dr -\int_t^T \tilde{Z}_rdW_r+\tilde{{K^l}}^o_T- \tilde{{K^u}}^o_T -(\tilde{{K^l}}^o_t- \tilde{{K^u}}^o_t).\enq
\end{Proposition}
\2
{\bf Proof. }
To develop this iteration, we suppose that $ (Y, Z, K^u, K^l) \in \mathcal{D}$, and $L \leq Y \leq U$, after the first iteration we obtain $ (\tilde{Y}, \tilde{Z}, \tilde{K}^u, \tilde{K}^l) \in \mathcal{D}$, and according to (\ref{KtSlaby})-(\ref{x}), we define
\begin{eqnarray} \label{Ktilde}
        \tilde{K}_t := \tilde{K}^l-\tilde{K}^u &=&  I_{\{T^U < T_L\}}I_{[T^U, \infty)}(t)\left[H_{L,U}(x_{T-t})-H_{L,U}(x_T) \right] \nonumber \\
        && + I_{\{T_L < T^U\}}I_{[T_L, \infty)}(t)\left[J_{L,U}(x_{T-t})-J_{L,U}(x_T) \right],
    \end{eqnarray}
    where
    \begin{equation}
    x_t = \xi  +\int_{T-t}^Tf(r, Y_r, Z_r, {K^l}^\flat_r-{K^u}^\flat_r)dr - \int_{T-t}^T Z_rdW_r.
    \end{equation}
Here $K^l$ and $K^u$ is replaced by the optional projections  ${K^l}^\flat$ and ${K^u}^\flat$ respectively.\\
To define $\tilde{Y}$, we first consider $\hat{Y}$ such that
\begin{eqnarray} \label{Yhat}
\hat{Y}_t=\xi + \int_t^T f(r, Y_r, Z_r, {K^l}^\flat_r-{K^u}^\flat_r)dr  + \tilde{{K^l}}_T- \tilde{{K^u}}_T - (\tilde{{K^l}}_t- \tilde{{K^u}}_t) - \int_t^TZ_rdW_r,
\end{eqnarray}
%where $\tilde{M}_t = \mathbb{E} \left( \xi + \tilde{{K^l}}^\flat_T - \tilde{{K^u}}^\flat_T+ \int_0^T f(r, Y_r, Z_r, {K^l}^\flat_r-{K^u}^\flat_r)dr / \mathcal{F}_t\right) - \tilde{N}_t$ is
%the martingale part such that $\tilde{N}_t = \tilde{{K^l}}^\flat_t - \tilde{{K^l}}^o_t - (\tilde{{K^u}}^\flat_t - \tilde{{K^u}}^o_t) $.
However, $\hat{Y}$ is not necessary adapted, the reason for which we consider its optional projection $\hat{Y}^\flat:=\tilde{Y}$ in the sense of Definition \ref{DefOP}, then
\beqs
\tilde{Y}_t &=& \E \left\{\xi + \int_t^T f(r, Y_r, Z_r, {K^l}^\flat_r-{K^u}^\flat_r)dr +\tilde{{K^l}}_T- \tilde{{K^u}}_T -(\tilde{{K^l}}_t- \tilde{{K^u}}_t) \Big/\Fc_t\right\} \\
            &=& M_t  -\int_0^t f(r, Y_r, Z_r, {K^l}^\flat_r-{K^u}^\flat_r)dr - (\tilde{{K^l}}^o_t- \tilde{{K^u}}^o_t ),
\enqs
where $M$ is a continuous martingale given by
\beqs M_t&=&\E \left\{\xi + \int_0^T f(r, Y_r, Z_r, {K^l}^\flat_r-{K^u}^\flat_r)dr +\tilde{{K^l}}_T- \tilde{{K^u}}_T  \Big/\Fc_t\right\}\\
                && -(\tilde{{K^l}}^\flat_t- \tilde{{K^l}}^o_t) + ( \tilde{{K^u}}^\flat_t- \tilde{{  K^u}}^o_t )\enqs
By martingale representation theorem, there exist a predictable process $\tilde{Z}$ such that
\beq  \tilde{Y}_t = \xi + \int_t^T f(r, Y_r, Z_r, {K^l}^\flat_r-{K^u}^\flat_r)dr -\int_t^T \tilde{Z}_rdW_r+\tilde{{K^l}}^o_T- \tilde{{K^u}}^o_T -(\tilde{{K^l}}^o_t- \tilde{{K^u}}^o_t).\enq
From Lipschitz property of $f$, and Proposition \ref{PropK} we have $(\tilde{Y}, \tilde{Z}, \tilde{K}^l,\tilde{K}^u )\in \Dc $. Moreover, the mapping $Y \rightarrow \tilde{Y}$ preserves the constraint $L \leq \tilde{Y} \leq U $. \ep
\begin{Remark} \label{RemKtild}
The process $\tilde{K}^l$ and  $\tilde{K}^u$ increase only on set $\{\hat{Y}=L\}$ and $\{\hat{Y}=U\}$ respectively.
\end{Remark}
\2
%%%%%%%%% Theorem   %%%%%%%%%%%%
We are going now to prove existence of the solution.
\begin{Theorem}
Assume that $(Y,Z,K^l,K^u)$ is a fixed point of $ \phi$, then $ (Y, Z, K^l,K^u)$ is a solution of the  reflected BSDE with resistance in the sense of Definition \ref{DefDRBSDER}. Moreover, the processes $K^l$ and $K^u$ are adapted.
\end{Theorem}
{\bf Proof. }
Since we suppose that $(Y,Z,K^l,K^u)$ is a fixed point of $\phi$, then $\phi(Y, Z, K^l, K^u) = (Y, Z, K^l, K^u)$ and according to Proposition \ref{phi} we have $ $:
\beq  Y_t = \xi + \int_t^T f(r, Y_r, Z_r, {K^l}^\flat_r-{K^u}^\flat_r)dr -\int_t^T Z_rdW_r+{K^l}^o_T- {K^u}^o_T -({K^l}^o_t- {K^u}^o_t)\enq
Then by (\ref{Ktilde}), it follows that
\begin{eqnarray}
        K^l_t - K^u_t &=&  I_{\{T^U < T_L\}}I_{[T^U, \infty)}\left[H_{L,U}(x_{T-t})-H_{L,U}(x_T) \right] \nonumber \\
        && + I_{\{T_L < T^U\}}I_{[T_L, \infty)}\left[J_{L,U}(x_{T-t})-J_{L,U}(x_T) \right],
    \end{eqnarray}
with \begin{equation}\label{xt}
x_t = \xi  +\int_{T-t}^Tf(r, Y_r, Z_r, {K^l}^\flat_r-{K^u}^\flat_r)dr - \int_{T-t}^T Z_rdW_r.
\end{equation}
By uniqueness of Skorohod equation, we have ${K^l}^o-{K^u}^o = K^l-K^u $, and since ${K^l}^o$ and ${K^l}^o$ are optional, it follows the adaptness of the process $K^l-K^u$, so that $K^l-K^u={K^l}^\flat-{K^u}^\flat$.
\ep  \vspace{5mm} \\
%%%%%%%%% T h e o r e m  %%%%%%%%%%%%
In the following, we prove the main result in this paper, which state uniqueness of the solution in the sense of Theorem 2.3 in Peng and Xu \cite{PX05}.
\begin{Theorem} \label{TheoremPointfix}
Assume Lipschitz continuity of $f$. There exists a unique quadriple $(Y, Z, K^l, K^u) $ in the space $\mathcal{D}$, solution to doubly reflected BSDE with resistance in the sense of Definition \ref{DefDRBSDER}.  The solution is unique in the following sense: if $(Y', Z', {K^l}', {K^u}) $ is another solution, then $Y\equiv Y' $, $ Z\equiv Z'$, and $ {K^l} - {K^u} \equiv {K^l}' -{K^u}'$, $\forall t\in [0, T]$.
\end{Theorem}
\s\\
{\bf Proof. } In order to simplify the notation, we not $\Kh^j =K^j-K^{'j}$ for $j=l,u$.
\2 \\
Let the space $\mathcal{D}:= \mathcal{S}^2 \times H^{2,d} \times \mathcal{A}^2 \times \mathcal{A}^2  $ be endowed with the norm:
\begin{eqnarray*}
\| (Y,Z,K^l,K^u) - (Y', Z', K^{'l},K^{'u}) \|^2_{\alpha, \beta} &:=&\| Y-Y'\|^2_{\alpha} + \| Z-Z'\|^2_{\alpha} + \beta \| \Kh^l - \Kh^{'u}\|^2_{\infty},
\end{eqnarray*}
where $\alpha \geq 0$, $\beta \geq 0$ and
\begin{eqnarray*}
\|\Kh^l - \Kh^{'u} \|^2_{\infty} &:=& \sup_{0 \leq t \leq T} \mathbb{E} |\Kh^l - \Kh^{'u} |^2 \\
\|Y \|^2_{\alpha} &:=& \int_0^T e^{\alpha r}\mathbb{E} |Y_r|^2 dr.
\end{eqnarray*}
Let the mapping $\phi$ defined in proposition \ref{phi}, be such that $\phi(Y, Z, K^l, K^u) =(\tilde{Y}, \tilde{Z}, \tilde{K}^l, \tilde{K}^u) $ and $\phi(Y', Z', K^{'l}, K^{'u}) =(\tilde{Y}', \tilde{Z}', \tilde{K}^{'l},\tilde{K}^{'u}) $.
From now, the proof will be divided into three steps.
\2 \\
{\bf Step 1.} We show that:
\begin{eqnarray}\label{step1}
 \|\tilde{Y} - \tilde{Y}'\|^2_\alpha + \|\tilde{Z}-\tilde{Z}'\|^2_\alpha \leq  \frac{2L_1}{\gamma_1}\left(\|Y - Y'\|^2_\alpha+ \|Z - Z'\|^2_\alpha\right) + \frac{2L_2}{\gamma_2}\|\Kh^\flat - \Kh^{'\flat}\|^2_\alpha.
\end{eqnarray}
We applied Itô's formula to $e^{\alpha t}(\tilde{Y}_t - \tilde{Y}'_t)^2 $, and taking its expectation we have,
\begin{eqnarray}
\E e^{\alpha t}(\tilde{Y}_t - \tilde{Y}'_t)^2 &=& -\alpha \E\int_t^T e^{\alpha r}(\tilde{Y}_r - \tilde{Y}'_r)^2 dr  - \E\int_t^Te^{\alpha r}|\tilde{Z}_r - \tilde{Z}'_r|^2dr\nonumber \\
        &&+ 2\E\int_t^Te^{\alpha r} (\tilde{Y}_r - \tilde{Y}'_r)d(\tilde{K}^{o}_r - \tilde{K}^{'o}_r ) \nonumber \\
        &&+ 2\E\int_t^Te^{\alpha r} (\tilde{Y}_r - \tilde{Y}'_r)\big[f(r, Y_r, Z_r, K^\flat_r) - f(r, Y'_r, Z'_r, K^{'\flat}))\big]dr ,\label{Yhat2}
\end{eqnarray}
where $ \K^o = \tilde{K^l}^o-\tilde{K^u}^o$ and $ \K^\flat = \tilde{K^l}^\flat-\tilde{K^u}^\flat$. Now, since $\tilde{Y}$ and $\tilde{Y}'$ are optional, then by Definition \ref{DefDOP}
\beqs \E\int_t^Te^{\alpha r} (\tilde{Y}_r - \tilde{Y}'_r)d(\tilde{K}^{o}_r - \tilde{K}^{'o}_r )&=&\E\int_t^Te^{\alpha r} (\tilde{Y}_r - \tilde{Y}'_r)d(\tilde{K}_r - \tilde{K}'_r ) \\
&=& \E\int_t^Te^{\alpha r} (\Y_r - \Y'_r)d(\K^l_r- \K^{l'}_r-\K^u_r +\K^{u'}_r ). \\
\enqs
Observe that
\beqs
\E\int_t^Te^{\alpha r} (\Y_r - \Y'_r)d(\K^l_r- \K^{l'}_r)
&=&  \E\int_t^Te^{\alpha r} (\Y_r - L_r)d\K^l_r - \E\int_t^Te^{\alpha r} (\Y_r - L_r)d\K^{l'}_r\\
&&- \E\int_t^Te^{\alpha r} (\Y'_r - L_r)d\K^l_r + \E\int_t^Te^{\alpha r} (\Y'_r - L_r)d\K^{l'}_r\\
&\leq&\E\int_t^Te^{\alpha r} (\Y_r - L_r)d\K^l_r+ \E\int_t^Te^{\alpha r} (\Y'_r - L_r)d\K^{l'}_r,
\enqs
and similarly
\beqs
-\E\int_t^Te^{\alpha r} (\Y_r - \Y'_r)d(\K^u_r- \K^{u'}_r) \\
&=& - \E\int_t^Te^{\alpha r} (\Y_r - U_r)d\K^u_r + \E\int_t^Te^{\alpha r} (\Y_r - U_r)d\K^{u'}_r\\
&&+ \E\int_t^Te^{\alpha r} (\Y'_r - U_r)d\K^u_r - \E\int_t^Te^{\alpha r} (\Y'_r - U_r)d\K^{u'}_r\\
&\leq&-\E\int_t^Te^{\alpha r} (\Y_r - U_r)d\K^u_r- \E\int_t^Te^{\alpha r} (\Y'_r - U_r)d\K^{u'}_r.
\enqs
Since $\K^l$ is increasing and $\tilde{Y}$ is the optional projection of $\hat{Y}$, then by Theorem 4.16 in \cite{N06}, and Remark \ref{RemKtild}
\beqs \E\int_t^Te^{\alpha r} (\Y_r - L_r)d\K^l_r &=& \E\int_t^Te^{\alpha r} (\hat{Y}_r - L_r)d\K^l_r\\
&=&0.
\enqs
Similarly, we show that  $ \E\int_t^Te^{\alpha r} (\Y'_r - L_r)d\K^{l'}_r=\E\int_t^Te^{\alpha r} (\Y_r - U_r)d\K^u_r= \E\int_t^Te^{\alpha r} (\Y'_r - U_r)d\K^{u'}_r=0$.
\2 \\
Plugging this in to (\ref{Yhat2}), and using the Lipschitz continuity of $f$,
\beqs\mathbb{E}\left(e^{\alpha t}(\tilde{Y}_t - \tilde{Y}'_t)^2 \right) &\leq& -\alpha \int_t^T\mathbb{E} \left(e^{\alpha r}(\tilde{Y}_r - \tilde{Y}'_r)^2\right) dr -\mathbb{E}\int_t^Te^{\alpha
        r}|\tilde{Z}_r-\tilde{Z}'_r|dr\nonumber \\
        &&+ 2L_1\int_t^Te^{\alpha r} \mathbb{E}\left(|\tilde{Y}_r - \tilde{Y}'_r|(|Y_r - Y'_r|+ |Z_r - Z'_r|)\right)dr\nonumber \\
        &&+ 2L_2\int_t^Te^{\alpha r} \mathbb{E}\left(|\tilde{Y}_r - \tilde{Y}'_r||K^{\flat}_r - K^{'\flat}_r|\right)dr.\nonumber
\enqs
Using the fact that $ ab \leq \gamma a^2 + \frac{1}{\gamma}b^2$ for some positive constant $\gamma$:
\begin{eqnarray}\label{step1'}
\mathbb{E}\left(e^{\alpha t}(\tilde{Y}_t - \tilde{Y}'_t)^2 \right) &\leq& -(\alpha -  \gamma_1L_1 -\gamma_2L_2   )\int_t^T\mathbb{E} \left(e^{\alpha
r}(\tilde{Y}_r - \tilde{Y}'_r)^2\right) dr \nonumber \\
        &&+  \frac{2L_1}{\gamma_1}\int_t^Te^{\alpha r} \mathbb{E}\left(|Y_r - Y'_r|^2+ |Z_r - Z'_r|^2\right)dr\nonumber \\
        &&+ \frac{2L_2}{\gamma_2}\int_t^Te^{\alpha r} \mathbb{E}|K^{\flat}_r - K^{'\flat}_r|^2 dr - \mathbb{E}\int_t^Te^{\alpha r}|\tilde{Z}_r-\tilde{Z}'_r|^2dr.
\end{eqnarray}
\vspace{2mm}
The result follows for $\alpha = \gamma_1 L_1 +\gamma_2 L_2$ and $t=0$.\2 \\
\noindent {\bf Step 2.}: We show that
\begin{eqnarray}\label{step2}
\|\tilde{K}_t - \tilde{K}'_t \|^2_\infty      &\leq& 45(TL_1^2+C_b)\left(\|Y-Y'\|^2_0 + \|Z-Z'\|^2_0\right) +45TL^2_2\|K^\flat-K^{'\flat}\|^2_\infty.
\end{eqnarray}
where $C_b$ is the constant appearing in Burkholder-Davis-Gundy inequality.
\2 \\
From Skorohod equation (\ref{SkorEq}), equation (\ref{Yhat}) and the fact that $K_t = \Xi_{L,U}(x_{T-t}) - \Xi_{L,U}(x_T) $, we have:
\beqs
\tilde{K}_t &=& -\int_0^t f(r, Y_r, Z_r, K^\flat_r)dr + \int_0^t Z_rdW_r + \ \Gamma_{L,U}(x_T) - \Gamma_{L,U}(x_{T-t}) \\
\tilde{K}'_t &=& -\int_0^t f(r, Y'_r, Z'_r, K^{'\flat}_r)dr + \int_0^t Z'_rdW_r  + \Gamma_{L,U}(x'_T) - \Gamma_{L,U}(x'_{T-t}),
\enqs
where
\begin{eqnarray*}
x_t &=& \xi  +\int_{T-t}^Tf(r, Y_r, Z_r, K^\flat_r)dr - \int_{T-t}^T Z_rdW_r, \\
x'_t &=& \xi  +\int_{T-t}^Tf(r, Y'_r, Z'_r, K^{'\flat}_r)dr - \int_{T-t}^T Z'_rdW_r.
\end{eqnarray*}
It follows
\begin{eqnarray*}
|\tilde{K}_t-\tilde{K}'_t|&\leq& \left|\int_0^t\left[f(r, Y_r, Z_r, K^{\flat}_r) - f(r, Y'_r, Z'_r, K^{'\flat}_r)dr \right]dr \right| \\
     && + | \Gamma_{L,U}(x_T) - \Gamma_{L,U}(x'_T)| + | \Gamma_{L,U}(x_{T-t}) - \Gamma_{L,U}(x'_{T-t})| \\
     &&+ \left|\int_0^t \left[Z_r-Z'_r\right]dW_r\right|.
     \end{eqnarray*}
Then
\begin{eqnarray} \label{K1}
|\tilde{K}_t-\tilde{K}'_t|^2&\leq& 3T\int_0^T\left|f(r, Y_r, Z_r, K^\flat_r) - f(r, Y'_r, Z'_r, K^{'\flat}_r)dr \right|^2dr \nonumber \\
     && + 6\| \Gamma_{L,U}(x) - \Gamma_{L,U}(x')\|^2 + 3\left|\int_0^T \left|Z_r-Z'_r\right|dW_r\right|^2.
\end{eqnarray}
From other side, we have by Theorem \ref{SlabyLipsh}:
\begin{eqnarray*}
\| \Gamma_{L,U}(x) - \Gamma_{L,U}(x')\|^2 &\leq& \left|\sup_{0 \leq t \leq T} | x_t - x'_t| \right|^2 \\
                                        &\leq& 2T\int_0^T\left|f(r, Y_r, Z_r, K^\flat_r) - f(r, Y'_r, Z'_r, K^{'\flat}_r)\right|^2dr \\
                                        &&+ 2\left|\int_0^T \left|Z_r-Z'_r\right|dW_r\right|^2.
\end{eqnarray*}
Plugging this in (\ref{K1}), and taking expectation in both hand side, we have:
\begin{eqnarray}
\mathbb{E}|\tilde{K}_t-\tilde{K}'_t|^2&\leq& 15T\mathbb{E}\int_0^T\left|f(r, Y_r, Z_r, K^\flat_r) - f(r, Y'_r, Z'_r, K^{'\flat}_r)dr \right|^2dr \nonumber \\
     && + 15\mathbb{E}\left|\int_0^T \left|Z_r-Z'_r\right|dW_r\right|^2.
\end{eqnarray}
Taking the supremum on $t$ over $[0, T]$, using Lipschitz continuity of $f$,  Jensen inequality together with Burkhoder-Davis-Gundy inequality:
\begin{eqnarray*}
\|\tilde{K}-\tilde{K}'\|^2_\infty &\leq& 15T\mathbb{E}\int_0^T\left|f(r, Y_r, Z_r, K^\flat_r) - f(r, Y'_r, Z'_r, K^{'\flat}_r)dr \right|^2dr \nonumber \\
     && + 15\mathbb{E}\left[\left|\int_0^T \left|Z_r-Z'_r\right|dW_r\right|^2\right]    \\
     &\leq& 45(TL_1^2+C_b)\mathbb{E}\int_0^T\left(|Y_r-Y'_r|^2 + |Z_r-Z'_r|^2\right)dr \nonumber \\
     && + 45TL^2_2\mathbb{E}\int_0^T|K^\flat_r-K^{'\flat}_r|^2 dr \\
    &\leq& 45(TL_1^2+C_b)\left(\|Y-Y'\|^2_0 + \|Z-Z'\|^2_0\right) +45TL^2_2\|K^\flat-K^{'\flat}\|^2_\infty.
\end{eqnarray*}
\noindent {\bf Step 3.} \\
Since  we have :
\begin{eqnarray}
\|K^\flat-K^{'\flat}\|^2_\alpha &\leq& \frac{e^{\alpha T } -1}{\alpha}\|K-K'\|^2_\infty.
\end{eqnarray}
Combining (\ref{step1}) and (\ref{step2}), leads to:
\begin{eqnarray}
  \|\tilde{Y} &-& \tilde{Y}'\|^2_\alpha + \|\tilde{Z}-\tilde{Z}'\|^2_\alpha + \beta \|\tilde{K} - \tilde{K}' \|^2_\infty  \nonumber \\
    &\leq& \left(\frac{2L_1}{\gamma_1}+45\beta(TL_1^2+C_b) \right)\left(\|Y-Y'\|_0 + \|Z-Z'\|_0\right) \nonumber \\
    &+& \left(\frac{2L_2(e^{\alpha T}-1)}{\alpha \beta \gamma_2 }+ 45TL^2_2\right)\beta\|K-K'\|_\infty.
\end{eqnarray}
We set $\gamma_1=2L_1^{-1} $, $\gamma_2=2L_2^{-1} $, $\alpha =5$ and we choose $L_1$, $L_2$ small enough such that
\beqs L_1 < \sqrt{\frac{\frac{1}{2}-C_b}{45T\beta}}, \s\s L_2 < \sqrt{\frac{1}{2}\times\frac{5\beta}{e^{5T}+225\beta-1}},\enqs
Then
\beqs\| (\tilde{Y},\tilde{Z},\tilde{K}^l,\tilde{K}^u) - (\tilde{Y}', \tilde{Z}', \tilde{K'}^l,\tilde{K'}^u) \|^2_{\alpha, \beta}  \leq \frac{1}{2} \| (Y,Z,K^l,K^u) - (Y', Z', K^{'l},K^{'u}) \|^2_{\alpha, \beta}. \enqs
So, the mapping is a contraction, and there is a fixed point $(Y, Z,K^l,K^u)$, which is the solution.
\ep \\

%%%%%%%%%%%%%%%%%%%%%%%%%%%%%%%%%%%%%%%%%%%%%%%%%%%%%%%%%%%%%%%%%%%%%%%%%%%%%%%%%%%%%%%%%%%%%
%%%%%%%%%%%%%%%%%%%%%% C o n t i n u o u s   d e p e n d a n c e                  %%%%%%%%%%%%
%%%%%%%%%%%%%%%%%%%%%%%%%%%%%%%%%%%%%%%%%%%%%%%%%%%%%%%%%%%%%%%%%%%%%%%%%%%%%%%%%%%%%%%%%%%%%
%
\section{Continuous dependance }
Our formulation of  doubly reflected BSDE with resistance permits us to derive the following continuous dependence thorem. 
\begin{Proposition}
Under the same assumptions in Theorem \ref{TheoremPointfix}. Let $(Y^i, Z^i, K^{il},K^{iu})$  with $i=1, 2$, be  solution of the following DRBSDE:
\begin{eqnarray*}
Y^i_t = \xi^i  +\int_t^Tf(r, Y^i_r, Z^i_r,K^i_r)dr - \int_t^T Z^i_rdW_r + K^{il}_T - K^{il}_t - (K^{iu}_T - K^{iu}_t)
\end{eqnarray*}
with two obstacles $L$ and $U$, where $K^i:= K^{il} - K^{iu}$, in the sense of Definition \ref{DefDRBSDER}.  Then we have
\begin{eqnarray*}
&\mathbb{E}&\left[\sup_{0\leq t \leq T}|Y^1_t - Y^2_t|^2 \right] + \mathbb{E} \left[\sup_{0\leq t \leq T}\Big|{K^1}^l_t - {K^2}^l_t - ({K^1}^u_t - {K^2}^u_t)\Big|^2 \right] \\
&& +\mathbb{E}\left[\int_0^T|Z^1_r - Z^2_r|^2dr \right]    \\
 &\leq&C\mathbb{E}\left[|\xi^1 - \xi^2|^2\right].
\end{eqnarray*}
The constant $C$ depends only on $L_1$, $L_2$ and T.
\end{Proposition}

{\bf Proof. }
We set $\hat{Y} = Y^1 - Y^2$, $\hat{Z} = Z^1 - Z^2$, $\hat{K}^l = {K^1}^l - {K^2}^l $,  $\hat{K}^u={K^1}^u - {K^2}^u$, $\hat{K}=K^1 - K^2 $$\hat{\xi} = \xi^1 - \xi^2$, $\hat{f}_r= f(r,Y^1,Z^1,K^1) - f(r,Y^2,Z^2,K^2)$, we have
\begin{eqnarray}\label{hat{Y}}
\hat{Y}_t = \hat{\xi}  +\int_t^T\hat{f}_rdr - \int_t^T \hat{Z}_rdW_r + \hat{K}^l_T - \hat{K}^l_t - (\hat{K}^u_T - \hat{K}^u_t)
\end{eqnarray}
Apply Itô Formula to $|\hat{Y}_t|^2$, then
\begin{eqnarray} \label{ItoY2}
|\hat{Y}_t|^2 + \int_t^T|\hat{Z}_r|^2dr  &=& |\hat{\xi}|^2  +2\int_t^T \hat{Y}_r \hat{f}_rdr + 2\int_t^T\hat{Y}_rd\hat{K}^l_r \nonumber \\
                            &-& 2\int_t^T\hat{Y}_rd\hat{K}^u_r - 2\int_t^T \hat{Y}_r \hat{Z}_rdW_r
\end{eqnarray}
We first observe that

\begin{eqnarray}
\int_t^T(Y^1_r - L_s)dK^{1l}_r=\int_t^T(Y^2_r - L_r)dK^{2l}_r=0 \nonumber \\
\int_t^T(Y^1_r - U_r)dK^{1u}_r=\int_t^T(Y^2_r - U_r)dK^{2u}_r=0
\end{eqnarray}
Thus
\begin{eqnarray}\label{dKl}
\int_t^T\hat{Y}_rd\hat{K}^l_r &=& \int_t^T (Y^1_r - L_r) dK^{1l}_r + \int_t^T (L_r - Y^2_r) dK^{1l}_r \nonumber \\
 && + \int_t^T (Y^2_r - L_r) dK^{2l}_r + \int_t^T (L_r - Y^1_r) dK^{2l}_r\nonumber \\
 &\leq& 0
\end{eqnarray}
And
\begin{eqnarray}\label{dKu}
\int_t^T\hat{Y}_rd\hat{K}^u_r &=& \int_t^T (Y^1_r - U_r) dK^{1u}_r + \int_t^T (U_r - Y^2_r) dK^{1u}_r \nonumber \\
 &+& \int_t^T (Y^2_r - U_r) dK^{2u}_r + \int_t^T (U_r - Y^1_r) dK^{2u}_r\nonumber \\
 &\geq& 0
\end{eqnarray}
Applying this to the equation (\ref{ItoY2}) we obtain:
\begin{eqnarray}
|\hat{Y}_t|^2 + \int_t^T|\hat{Z}_r|^2dr  \leq  |\hat{\xi}|^2  +2\int_t^T \hat{Y}_r \hat{f}_rdr  - 2\int_t^T \hat{Y}_r \hat{Z}_rdW_r
\end{eqnarray}
By Lipschitz condition of $f$ we have:
\begin{eqnarray*}
|\hat{Y}_t|^2 + \int_t^T|\hat{Z}_r|^2dr  &\leq&  |\hat{\xi}|^2  +2\int_t^T \hat{Y}_r (L_1(|\hat{Y}_r| + |\hat{Z}_r|) + L_2|\hat{K}_r|)dr  \\
                                        &-& 2\int_t^T \hat{Y}_r \hat{Z}_rdW_r \\
                                        &\leq& |\hat{\xi}|^2  +(2L_1 + \alpha L^2_1 + \beta)\int_t^T |\hat{Y}_r|^2dr + \frac{1}{\alpha} \int_t^T |\hat{Z}_r|^2dr \\
                                        &+& \frac{L^2_2}{\beta}\int_t^T |\hat{K}_r|^2dr - 2\int_t^T \hat{Y}_r \hat{Z}_rdW_r
\end{eqnarray*}
Set $\alpha =2$ we have
\begin{eqnarray*}
\mathbb{E}\left[|\hat{Y}_t|^2\right]  &\leq&  \mathbb{E}\left[|\hat{\xi}|^2 \right]   +(2L_1 + 2 L^2_1 + \beta)\mathbb{E}\int_t^T |\hat{Y}_r|^2dr + \frac{TL^2_2}{\beta}\mathbb{E}\Big[\sup_{0 \leq t \leq T}\hat{K}^2_r\Big]
\end{eqnarray*}
By Gronwall Lemma we have:
\begin{eqnarray*}
\mathbb{E}\left[|\hat{Y}_t|^2\right]  &\leq&   e^{2L_1 + 2 L^2_1 + \beta}\mathbb{E}\left[|\hat{\xi}|^2 + \frac{TL^2_2}{\beta}\sup_{0 \leq t \leq T}\hat{K}^2_r\right]
\end{eqnarray*}
it follows that
\begin{eqnarray*}\label{YtZt}
\mathbb{E}\left[|\hat{Y}_t|^2 + \int_t^T|\hat{Z}_r|^2dr  \right]  &\leq&   e^{2L_1 + 2 L^2_1 + \beta}\mathbb{E}\left[|\hat{\xi}|^2 + \frac{TL^2_2}{\beta}\sup_{0 \leq t \leq T}\hat{K}^2_r \right]
\end{eqnarray*}

Using Burkholder-Davis-Gundy inequality, we have
\begin{eqnarray} \label{supY}
\mathbb{E}\left[\sup_{0\leq t \leq T}|\hat{Y}_t|^2 + \int_0^T|\hat{Z}_r|^2dr  \right]  &\leq&   e^{2L_1 + 2 L^2_1 + \beta}\mathbb{E}\left[|\hat{\xi}|^2 + \frac{TL^2_2}{\beta}\sup_{0 \leq t \leq T}\hat{K}^2_r \right] \nonumber \\
\end{eqnarray}
From other side, we have by equation (\ref{hat{Y}}):
\begin{eqnarray*}
\hat{K}_t = \hat{Y}_0 - \hat{Y}_t  - \int_0^t\hat{f}_rdr + \int_0^t \hat{Z}_rdW_r
\end{eqnarray*}
Using again Burkholder-Davis-Gundy inequality,, Jensen inequality and  (\ref{supY}) :
\beqs
\sup_{0 \leq t \leq T}|\hat{K}_t|^2   &\leq&  C\left( \sup_{0 \leq t \leq T}|\hat{Y}_t|^2 + T\int_0^T \Big[L_1^2\big(|\hat{Y}_r|^2+|\hat{Z}_r|^2\big)+L_2^2|\hat{K}_r|^2\Big]dr \right. \\
                && \s +\left. \sup_{0 \leq t \leq T}\Big(\int_0^t\hat{Z}_rdW_r \Big)^2\right)
\enqs
Set $\beta=1$, leads to
\beqs
\mathbb{E}\sup_{0 \leq t \leq T}|\hat{K}_t|^2 &\leq&  Ce^{2L_1 + 2 L^2_1 + \beta}(2 + TL_1^2 + T^2L_1^2)\E[\hat{\xi}^2] \\
                    &&+ CTL_2^2\big(e^{2L_1 + 2 L^2_1 + \beta}(2 + TL_1^2 + T^2L_1^2)+1 \big)\mathbb{E}\sup_{0 \leq t \leq T}|\hat{K}_t|^2
\enqs
Choosing $L_1$ and $L_2$  small enough we have
\beqs
\mathbb{E}\sup_{0 \leq t \leq T}|\hat{K}_t|^2 &\leq&  Ce^{2L_1 + 2 L^2_1 + \beta}(2 + TL_1^2 + T^2L_1^2)\E[\hat{\xi}^2]
\enqs
Plugging this in (\ref{supY}), we get the required result.\ep

%===================================== REFERENCES  =====================================


\begin{thebibliography}{1}
\bibitem {AO11}     Aazizi S. and Y. Ouknine (2011). Portfolio-constrained Backward SDEs with Jump and related Integro-partial differential equation, preprint.

\bibitem {BK04}     Bank P. and El N. Karoui (2004). A stochastic representation theorem with applications to optimization and obstacle problems, Ann. Probab. 32, no. 1B, 1030-1067. MR 2044673 (2005a:60046).

\bibitem {B76}      Bismut J.M (1976). Théorie probabiliste du contrôle des diffusions, Mem. Amer. Math.Soc. 4, no. 167, xiii+130. MR 0453161 (56 11428).

\bibitem {BKR09}    Burdzy K., W. Kang and K. Ramanan (2009). The Skorokhod problem in a time-dependent interval, Stochastic Processes and Their Applications , vol. 119, no. 2, pp. 428-452.

\bibitem {CK96}     Cvitanic J. and I. Karatzas (1996). Backward SDEs with reflection and Dynkin games, Annals of Probability 24 (4), pp. 2024-2056.

\bibitem {CKS98}    Cvitanic J., I. Karatzas and M. Soner (1998). Backward stochastic differential equations with constraints on the gain-process", Annals of Probability, 26, 1522-1551.

\bibitem {EC78}     El Karoui N. and M. Chaleyat-Maurel (1978). Un problème de réflexion et ses applications au temps local et aux equations differentielles
    stochastiques sur $\mathbb{R}$. Cas continu. In Temps Locaux. Astérisque 52-53 117-144. Soc. Math. France, Paris.

\bibitem {EKPPQ97}  El Karoui N., C. Kapoudjian, E. Pardoux,  S. Peng  and M. C. Quenez (1997). Reflected solutions of backward SDE's, and related obstacle problems for PDE's, Ann. Probab. 25, no. 2, 702-737. MR 1434123 (98k:60096).

\bibitem {HH06}     Hamadène S. and M. Hassani (2005). BSDEs with two reflecting barriers: the general result, Probability Theory and Related Fields 132, (2005), pp.237-264.

\bibitem{HWY92}     He S. W., J. G. Wang and J. A. Yan (1992). Semimartingales and Stochastic Calculus, CRC, Press nd Science Press.

\bibitem {KLRS07}   Kruk L., J. Lehoczky, K. Ramanan  and S. Shreve (2007). An explicit formula for the Skorohod map on $[0,a]$, Annals of Probability, vol. 35, no. 5, pp. 1740-1768.

\bibitem {MW09}     Ma J. and Y. Wang (2009). On variant reflected backward SDEs, with applications, J.Appl. Math. Stoch. Anal., Art. ID 854768, 26. MR 2511615 (2010g:60139).

\bibitem {N06}      Nikeghbali A. (2006). An essay on the general theory of stochastic processes, Probability Surveys Vol. 3, 345-412, ISSN: 1549-5787, DOI: 10.1214/154957806000000104.

\bibitem {PP90}     Pardoux E. and Peng S. (1990). Adapted solution of a backward stochastic differential equation, Systems Control Lett. 14, no. 1, 55-61. MR 1037747 (91e:60171).

\bibitem {PX05}     Peng S. and M. Xu (2005). The smallest g-supermartingale and reflected BSDE with single and double $L^2$-obstacles, Ann. I. H. Poincaré PR 41, 605-630.
\bibitem{R06}       Ramanan K. (2006). Reflected diffusions defined via the extended Skorokhod map, Electron. J. Probab. 11, 934-992.

\bibitem {S10}      Slaby M. (2010). An Explicit Representation of the Extended Skorohod Map with Two Time-Dependent Boundaries, Journal of Probability and Statistics Volume, Article ID 846320, 18 pages doi:10.1155/2010/846320.

\bibitem {QX11}     Qian Z. and M. Xu (2011). Skorohod Equation and Reflected Backward Stochastic Differential Equations, http://arxiv.org/abs/1103.2078.
\end{thebibliography}
\end{document}